\documentclass[11pt]{amsart}
\usepackage{amssymb,latexsym,amsmath,graphicx}

\theoremstyle{plain}
\newtheorem{theorem}{Theorem}
\numberwithin{theorem}{section}
\newtheorem{lemma}[theorem]{Lemma}

\theoremstyle{definition}
\newtheorem{definition}[theorem]{Definition}

\begin{document}

\title{Legendrian Vertical Circles in Small Seifert Spaces}

\author{Hao Wu}

\address{Department of mathematics\\
         MIT, Room 2-487\\
         77 Massachusetts Avenue\\
         Cambridge, MA 02139\\
         USA}

\email{haowu@math.mit.edu}

\maketitle

\begin{abstract}
We discuss the relations between the $e_0$ invariants of a small
Seifert space and the twisting numbers of Legendrian vertical
circles in it.
\end{abstract}

\section{Introduction}

A contact structure $\xi$ on an oriented $3$-manifold $M$ is a
nowhere integrable tangent plane distribution, i.e., near any
point of $M$, $\xi$ is defined locally by a $1$-form $\alpha$,
s.t., $\alpha\wedge d\alpha\neq0$. Note that the orientation of
$M$ given by $\alpha\wedge d\alpha$ depends only on $\xi$, not on
the choice of $\alpha$. $\xi$ is said to be positive if this
orientation agrees with the native orientation of $M$, and
negative if not. A contact structure $\xi$ is said to be
co-orientable if $\xi$ is defined globally by a $1$-form $\alpha$.
Clearly, an co-orientable contact structure is orientable as a
plane distribution, and a choice of $\alpha$ determines an
orientation of $\xi$. Unless otherwise specified, all contact
structures in this paper will be co-oriented and positive, i.e.,
with a prescribed defining form $\alpha$ such that $\alpha\wedge
d\alpha>0$. A curve in $M$ is said to be Legendrian if it is
tangent to $\xi$ everywhere. $\xi$ is said to be overtwisted if
there is an embedded disk $D$ in $M$ such that $\partial D$ is
Legendrian, but $D$ is transversal to $\xi$ along $\partial D$. A
contact structure that is not overtwisted is called tight.
Overtwisted contact structures appear to be very "soft". It is
proven by Eliashberg in \cite{E1} that two overtwisted contact
structures are isotopic \textit{iff} they are homotopic as tangent
plane distributions. Tight contact structures are more rigid.
Classifications of tight contact structures up to isotopy are only
known for very limited classes of $3$-manifolds. (See, e.g.,
\cite{E2}, \cite{E3}, \cite{EH}, \cite{Gh}, \cite{GS}, \cite{H1},
\cite{H2}, \cite{K}.)

A small Seifert space is a Seifert fibred space with $3$ singular
fibers over $S^2$. Any regular fiber $f$ in a small Seifert space
$M(\frac{q_1}{p_1},\frac{q_2}{p_2},\frac{q_3}{p_3})$ admits a
canonical framing given by pulling back an arc in the base $S^2$
containing the projection of $f$. An embedded circle in
$M(\frac{q_1}{p_1},\frac{q_2}{p_2},\frac{q_3}{p_3})$ is said to be
vertical if it is isotopic to a regular fiber. Any vertical circle
inherits a canonical framing from the canonical framing of regular
fibers. We call this framing $Fr$.

\begin{definition}
Let $\xi$ be a contact structure on a small Seifert
space $M(\frac{q_1}{p_1},\frac{q_2}{p_2},\frac{q_3}{p_3})$, and
$L$ a Legendrian vertical circle in $(M,\xi)$. The twisting number
$t(L)$ of $L$ is the twisting number of $\xi|_L$ along $L$
relative to the canonical framing $Fr$ of $L$.
\end{definition}

In \cite{CGH}, Colin, Giroux and Honda divided the tight contact
structures on a small Seifert space into two types: those for
which there exists a Legendrian vertical circle with twisting
number $0$, and those for which no Legendrian vertical circles
with twisting number $0$ exist. It is proven in \cite{CGH} that,
up to isotopy, the number of tight contact structures of the first
type is always finite, and, unless the small Seifert space is also
a torus bundle, the number of tight contact structures of the
second type is also finite. Their work gives in principle a method
to estimate roughly the upper bound of the number of tight contact
structures on a small Seifert space. In the present paper, we
demonstrate that most small Seifert spaces admit only one of the
two types of tight contact structures. To make our claim precise,
we need the following invariant. (See, e.g., \cite{Go}.)

\begin{definition}
For a small Seifert space
$M=M(\frac{q_1}{p_1},\frac{q_2}{p_2},\frac{q_3}{p_3})$, define
$e_0(M) = \lfloor\frac{q_1}{p_1}\rfloor +
\lfloor\frac{q_2}{p_2}\rfloor + \lfloor\frac{q_3}{p_3}\rfloor$,
where $\lfloor x \rfloor$ is the greatest integer not greater than
$x$.
\end{definition}

Clearly, $e_0(M)$ is an invariant of $M$, i.e., it does not depend
on the choice of the representatives
$(\frac{q_1}{p_1},\frac{q_2}{p_2},\frac{q_3}{p_3})$. Now we can
formulate our claim precisely in the following theorems.

\begin{theorem}\label{0}
Let $M=M(\frac{q_1}{p_1},\frac{q_2}{p_2},\frac{q_3}{p_3})$ be a
small Seifert space. If $e_0(M)\geq0$, then every tight contact
structure on $M$ admits a Legendrian vertical circle with twisting
number $0$.
\end{theorem}

\begin{theorem}\label{-2}
Let $M=M(\frac{q_1}{p_1},\frac{q_2}{p_2},\frac{q_3}{p_3})$ be
a small Seifert space. If $e_0(M)\leq-2$, then no tight contact
structures on $M$ admit Legendrian vertical circles with twisting
number $0$.
\end{theorem}

In particular, Theorem \ref{-2} means that, for any small Seifert
space $M=M(\frac{q_1}{p_1},\frac{q_2}{p_2},\frac{q_3}{p_3})$,
either $M$ or $-M$ does not admit tight contact structures for
which there exists a Legendrian vertical circle with twisting
number $0$, where $-M$ is $M$ with reversed orientation. This is
because that $e_0(M)+e_0(-M)=-3$, and, hence, one of $e_0(M)$ and
$e_0(-M)$ has to be less than or equal to $-2$.

It turns out that the case when $e_0(M)=-1$ is the most difficult.
Only very weak partial results are known. For example, in
\cite{GS}, Ghiggini and Sch\"{o}nenberger proved that, when
$r\leq\frac{1}{5}$, no tight contact structures on the small
Seifert space $M(r,\frac{1}{3},-\frac{1}{2})$ admit Legendrian
vertical circles with twisting number $0$.

We have the following results about the case when $e_0(M)=-1$.

\begin{theorem}\label{-1}
Let $M=M(\frac{q_1}{p_1},\frac{q_2}{p_2},\frac{q_3}{p_3})$ be a
small Seifert space such that $0<q_1<p_1$, $0<q_2<p_2$ and
$-p_3<q_3<0$.

(1) If $\frac{q_1}{p_1}+\frac{q_3}{p_3}\geq0$ or
$\frac{q_2}{p_2}+\frac{q_3}{p_3}\geq0$ or
$\frac{q_1}{p_1}+\frac{q_2}{p_2}\geq1$, then every tight contact
structure on $M$ admits a Legendrian vertical circle with twisting
number $0$.

(2) If $q_3=-1$, $\frac{q_1}{p_1}<\frac{1}{2p_3-1}$ and
$\frac{q_2}{p_2}<\frac{1}{2p_3}$, then no tight contact structures
on $M$ admit Legendrian vertical circles with twisting number $0$.

(3) If $q_1=q_2=1$ and $p_1,p_2>-2\lfloor\frac{p_3}{q_3}\rfloor$,
then no tight contact structures on $M$ admit Legendrian vertical
circles with twisting number $0$.
\end{theorem}

To understand the proofs in this paper, readers need to be
familiar with the techniques developed by Giroux in \cite{Gi1} and
Honda in \cite{H1}. For those who are not, there is a concise
introduction to these techniques in \cite{GS}.

\vspace{.4cm}

\section{The $e_0\geq0$ Case}

The $e_0\geq0$ case is the simplest case. Theorem \ref{0} is a
special case of Lemma \ref{existence} below, which also implies
part (1) of Theorem \ref{-1}.

In the rest of this paper, $\Sigma$ means a three hole sphere;
$\Sigma_0$ means a properly embedded three hole sphere in $\Sigma
\times S^1$ isotopic to $\Sigma \times \{\text{pt}\}$. Let
$-\partial\Sigma\times S^1=T_1+T_2+T_3$, where the "$-$" sign
means reversing the orientation. We identify $T_i$ to
$\bf{R}^2/\bf{Z}^2$ by identifying the corresponding component of
$-\partial\Sigma\times\{\text{pt}\}$ to $(1,0)^T$, and
$\{\text{pt}\}\times S^1$ to $(0,1)^T$. An embedded circle in
$T_i$ is called vertical if it is essential and has slope
$\infty$.

The following lemma is purely technical.

\begin{lemma}\label{cutting}
Let $\xi$ be a tight contact structure on $\Sigma \times S^1$.
Assume that each $T_i$ is convex with two dividing curves of slope
$s_i$. Then there exist collar neighborhoods $T_1\times I$ and
$T_2\times I$ of $T_1$ and $T_2$ and a properly embedded vertical
convex annulus $A$ in $(\Sigma \times S^1)\setminus (T_1\times I
\cup T_2\times I)$ connecting $T_1\times\{1\}$ to $T_2\times\{1\}$
with Legendrian boundary satisfying that following:
\begin{enumerate}
    \item $T_1\times I$ and $T_2\times I$ are mutually disjoint and
disjoint from $T_3$;
    \item for $i=1,2$, $T_i\times\{0\}=T_i$ and $T_i\times\{1\}$ is convex with two
dividing curves of slope $s'_i \leq s_i$;
    \item $A$ has no $\partial$-parallel dividing curves, and the
Legendrian boundary of $A$ intersects the dividing sets of
$T_1\times\{1\}$ and $T_2\times\{1\}$ efficiently.
\end{enumerate}
\end{lemma}

\begin{proof}
If both $s_1$ and $s_2$ are $\infty$, then we can isotope $T_1$
and $T_2$ slightly to make them to have vertical Legendrian
divides. Connect a Legendrian divide of $T_1$ to a Legendrian
divide of $T_2$ by a properly embedded vertical convex annlus $A$.
Then we are done. If $s_1=\infty$, but $s_2$ is finite, then we
make $T_1$ to have vertical Legendrian divides, and $T_2$ to have
vertical Legendrian rulings. Connect a Legendrian divide of $T_1$
to a Legendrian ruling of $T_2$ by a properly embedded vertical
convex annlus $B$. Then no dividing curves of $B$ intersects
$B\cap T_1$. And we can decrease $s_2$ to $\infty$ by isotoping
$T_2$ across the dividing curves of $B$ starting and ending on
$B\cap T_2$ through bypass adding. We can keep $T_2$ disjoint from
both $T_1$ and $T_3$ through out this isotopy since bypass adding
can be done in a small neighborhood of the bypass and the original
surface. The we are back to the case when $s_1$ and $s_2$ are both
$\infty$.

Assume $s_i=\frac{q_i}{p_i}$ is finite for $i=1,2$, where $p_i>0$.
First, we isotope $T_1$ and $T_2$ slightly so that they have
vertical Legendrian rulings. Note that the Legendrian rulings
always intersect dividing curves efficiently. Then connect a
Legendrian ruling of $T_1$ to a Legendrian ruling of $T_2$ by a
properly embedded vertical convex annulus $A$ in $\Sigma \times
S^1$. If $A$ has no $\partial$-parallel dividing curves, then we
are done. If $A$ has a $\partial$-parallel dividing curve, say on
the $T_1$ side, then, after possibly isotoping $A$ slightly, we
can assume there is a bypass of $T_1$ on $A$. Adding this bypass
to $T_1$, we get a convex torus $T'_1$ with two dividing curves in
$\Sigma \times S^1$ that co-bounds a collar neighborhood of $T_1$.
We can make $T'_1$ to have vertical Legendrian ruling. By Lemma
3.5 of \cite{H1}, we have that $T'_1$ has two dividing curves of
the slope $s'_1=\frac{q'_1}{p'_1}<s_1$, where $0\leq p'_1<p_1$.
Now we delete the thickened torus between $T_1$ and $T'_1$ from
$\Sigma \times S^1$, and repeat the procedure above. This whole
process will stop in less than $p_1+p_2$ steps, i.e, we can either
find the collar neighborhoods and the annulus with the required
properties, or force one of $s_1$ and $s_2$ to decrease to
$\infty$. But the lemma is proved in the latter case. This
finishes the proof.
\end{proof}

\begin{lemma}\label{existence}
Let $M=M(\frac{q_1}{p_1},\frac{q_2}{p_2},\frac{q_3}{p_3})$ be a
small Seifert space such that $\frac{q_1}{p_1},\frac{q_2}{p_2}>0$
and $\frac{q_1}{p_1}+\frac{q_3}{p_3}\geq0$. Then every tight
contact structure on $M$ admits a Legendrian vertical circle with
twisting number $0$.
\end{lemma}

\begin{proof}
For $i=1,2,3$, let $V_i=D^2\times S^1$, and identify $\partial
V_i$ with $\bf{R}^2/\bf{Z}^2$ by identifying a meridian $\partial
D^2 \times \{\text{pt}\}$ with $(1,0)^T$ and a longitude
$\{\text{pt}\}\times S^1$ with $(0,1)^T$. Choose $u_i, v_i \in
\bf{Z}$ such that $0<u_i<p_i$ and $p_iv_i+q_iu_i=1$ for $i=1,2,3$.
Define an orientation preserving diffeomorphism
$\varphi_i:\partial V_i\rightarrow T_i$ by
\[
\varphi_i ~=~
\left(%
\begin{array}{cc}
  p_i & u_i \\
  -q_i & v_i \\
\end{array}%
\right),
\]
where $T_i$ and the coordinates on it are defined above. Then
\[
M=M(\frac{q_1}{p_1},\frac{q_2}{p_2},\frac{q_3}{p_3})\cong(\Sigma
\times S^1)\cup_{(\varphi_1\cup\varphi_2\cup\varphi_3)}(V_1\cup
V_2\cup V_3).
\]
Let $\xi$ be a tight contact structure on $M$. We first isotope
$\xi$ to make each $V_i$ a standard neighborhood of a Legendrian
circle $L_i$ isotopic to the $\frac{q_i}{p_i}$-singular fiber with
twisting number $n_i<0$, i.e., $\partial V_i$ is convex with two
dividing curves each of which has slope $\frac{1}{n_i}$ when
measured in the coordinates of $\partial V_i$ given above. Let
$s_i$ be the slope of the dividing curves of $T_i=\partial V_i$
measured in the coordinates of $T_i$. Then we have that
\[
s_i ~=~ \frac{-n_iq_i+v_i}{n_ip_i+u_i} ~=~
-\frac{q_i}{p_i}+\frac{1}{p_i(n_ip_i+u_i)} ~<~ -\frac{q_i}{p_i}.
\]

By Lemma \ref{cutting}, we can thicken $V_1$ and $V_2$ to $V'_1$
and $V'_2$ such that
\begin{enumerate}
    \item $V'_1$, $V'_2$ and $V_3$ are pairwise disjoint;
    \item for $i=1,2$, $T'_i=\varphi_i(\partial V'_i)$
    is convex with two dividing curves of slope
    $s'_i=-\frac{q'_i}{p}\leq s_i$, where $p,q'_i>0$;
    \item there exists a properly embedded vertical convex annulus
    $A$ connecting $T'_1$ to $T'_2$ that has no $\partial$-parallel dividing curves,
    and the Legendrian boundary of $A$ intersects the dividing sets of
    these tori efficiently.
\end{enumerate}

If none of the dividing curves of $A$ is an arc connecting the two
components of $\partial A$, then, by the Legendrian Realization
Principle (\cite{Gi1}, \cite{H1}), we can isotope $A$ to make a
vertical circle $L$ on $A$ which is disjoint from the dividing
curves Legendrian. Note that $A$ gives the canonical framing of
$L$, and the twisting number of $\xi|_L$ relative to $TA|_L$ is
$0$ by Proposition 4.5 of \cite{K2}. So $t(L)=0$.

If there are dividing curves connecting the two components of
$\partial A$, then cut $M\setminus(V'_1\cup V'_2\cup V_3)$ open
along $A$. We get an embedded thickened torus $T_3\times I$ in $M$
such that $T_3\times\{0\}=T_3$, and $T_3\times\{1\}$ is convex
with two dividing curves of slope $s'_3=\frac{q'_1+q'_2-1}{p}$.
Note that
\[
s'_3 ~=~ \frac{q'_1+q'_2-1}{p} ~\geq~ \frac{q'_1}{p} ~\geq~ -s_1
~>~ \frac{q_1}{p_1} ~\geq~ -\frac{q_3}{p_3} ~>~ s_3.
\]
According to Theorem 4.16 of \cite{H1}, there exists a convex
torus $T$ in $T_3\times I$ parallel to $T_3$ with vertical
dividing curves. We can then isotope $T$ to make it in standard
form. Then a Legendrian divide of $T$ is a Legendrian vertical
circle with twisting number $0$.
\end{proof}

\noindent\textit{Proof of Theorem \ref{0} and Theorem
\ref{-1}(1).} If
$M=M(\frac{q_1}{p_1},\frac{q_2}{p_2},\frac{q_3}{p_3})$ satisfies
that $e_0(M)\geq0$, then we can assume that $\frac{q_i}{p_i}>0$
for $i=1,2,3$. It's then clear that
$\frac{q_1}{p_1},\frac{q_2}{p_2}>0$ and
$\frac{q_1}{p_1}+\frac{q_3}{p_3}>0$. Thus, Lemma \ref{existence}
implies Theorem \ref{0}.

Now we assume
$M=M(\frac{q_1}{p_1},\frac{q_2}{p_2},\frac{q_3}{p_3})$ is a small
Seifert space such that $0<q_1<p_1$, $0<q_2<p_2$ and $-p_3<q_3<0$.
If $\frac{q_1}{p_1}+\frac{q_3}{p_3}\geq0$ or
$\frac{q_2}{p_2}+\frac{q_3}{p_3}\geq0$, then Lemma \ref{existence}
applies directly. If $\frac{q_1}{p_1}+\frac{q_2}{p_2}\geq1$, we
apply Lemma \ref{existence} to
$M=M(\frac{q_1}{p_1},\frac{q_3}{p_3}+1,\frac{q_2}{p_2}-1)$. This
proves Theorem \ref{-1}(1). \hfill $\Box$

\vspace{.4cm}

\section{The $e_0\leq-2$ Case}

\begin{definition}
Let $\xi$ be a contact structure on $\Sigma \times S^1$. $\xi$ is
said to be inappropriate if $\xi$ is overtwisted, or there exists
an embedded $T^2 \times I$ with convex boundary and $I$-twisting
at least $\pi$ such that $T^2 \times \{0\}$ is isotopic to one of
the $T_i$'s. $\xi$ is called appropriate if it is not
inappropriate.
\end{definition}

\begin{lemma}\label{appropriate}
Let $M=M(\frac{q_1}{p_1},\frac{q_2}{p_2},\frac{q_3}{p_3})$ be a
small Seifert space, and $\xi$ a tight contact structures on $M$.
Suppose that $V_1$, $V_2$, $V_3$ are tubular neighborhoods of the
three singular fibers, and $\Sigma \times S^1=M\setminus(V_1 \cup
V_2\cup V_3)$. Then $\xi|_{\Sigma \times S^1}$ is appropriate.
\end{lemma}

\begin{proof}
Without loss of generality, we assume $\partial V_i$ is identified
with $T_i$ by the diffeomorphism $\varphi_i$. $\xi|_{\Sigma \times
S^1}$ is clearly tight. If it is inappropriate, then there exists
an embedded $T^2 \times I$ with convex boundary and $I$-twisting
at least $\pi$ such that $T^2 \times \{0\}$ is isotopic to one of
the $T_i$'s. Let's say $T^2 \times \{0\}$ is isotopic to $T_1$.
$T^2 \times I$ has $I$-twisting at least $\pi$ implies that, for
any rational slope $s$, there is a convex torus $T_0$ contained in
$T^2 \times I$ isotopic to $T_1$ that has dividing curves of slope
$s$. Specially, we let $m$ be a meridian of $\partial V_1$, and
$s$ the slope of $\varphi_1(m)$. Then the above fact means that we
can thicken $V_1$ so that $\partial V_1$ has dividing curves
isotopic to its meridians, which implies that the thickened $V_1$
is overtwisted. This contradicts the tightness of $\xi$. Thus,
$\xi|_{\Sigma \times S^1}$ must be appropriate.
\end{proof}

\begin{lemma}[\cite{EH}, Lemma 10]
Let $\xi$ be an appropriate contact structure on $\Sigma \times
S^1$ such that all three boundary components of $\Sigma \times
S^1$ are convex, and the dividing set of each of these boundary
components consists of two vertical circles. If $\Sigma_0$ is
convex with Legendrian boundary that intersects the dividing set
of $\partial\Sigma \times S^1$ efficiently, then the dividing set
of $\Sigma_0$ consists of three properly embedded arcs, each of
which connects a different pair of components of
$\partial\Sigma_0$.
\end{lemma}

The following lemma is a special case of Proposition 6.4 of
\cite{CGH}, which also appears in \cite{Gi2} and \cite{H2}.

\begin{lemma}[\cite{CGH}, \cite{Gi2}, \cite{H2}]\label{center}
Isotopy classes of tight contact structures on $\Sigma\times S^1$
such that all three boundary components of $\Sigma\times S^1$ are
convex, and the dividing set of each of these boundary components
consists of two vertical circles are in 1-1 correspondence with
isotopy classes of embedded multi-curves on $\Sigma$ with 2 fixed
end points on each component of $\partial\Sigma$ that have no
homotopically trivial components.

The correspondence here is induced by the following mapping: Given
a tight contact structure $\xi$ on $\Sigma\times S^1$, isotope
$\Sigma_0$ to make it convex with Legendrian boundary intersecting
the dividing set of $\partial\Sigma \times S^1$ efficiently. Then
map $\xi$ to the isotopy class of the dividing curves of
$\Sigma_0(\cong \Sigma$).
\end{lemma}

The following lemma from \cite{GS} plays a key role in the proof
of Theorem 1. For the convenience of readers, we give a detailed
proof here.

\begin{lemma}[\cite{GS}, Lemma 36]\label{GS36} Let $\xi$
be an appropriate contact structure on $\Sigma\times S^1$. Suppose
that $-\partial\Sigma\times S^1=T_1+T_2+T_3$ is convex and such
that each of $T_1$ and $T_2$ has vertical Legendrian rulings and
two dividing curves of slope $-\frac{1}{n}$, where $n\in{\bf
Z}^{>0}$, and $T_3$ has two vertical dividing curves. Let
$T_1\times I$ and $T_2\times I$ be collar neighborhoods of $T_1$
and $T_2$ that are mutually disjoint and disjoint from $T_3$, and
such that, for $i=1,2$, $T_i\times\{0\}=T_i$ and $T_i\times\{1\}$
is convex with dividing set consisting of two  vertical circles.
If $\xi|_{T_1\times I}$ and $\xi|_{T_2\times I}$ are both isotopic
to a given minimal twisting tight contact structure $\eta$ on
$T^2\times I$ relative to the boundary, then there exists a
properly embedded convex vertical annulus $A$ with no
$\partial$-parallel dividing curves, whose Legendrian boundary
$\partial A=(A\cap T_1)\cup (A\cap T_2)$ intersects the dividing
sets of $T_1$ and $T_2$ efficiently.
\end{lemma}

\begin{proof}
Let $\Sigma'\times S^1=(\Sigma\times S^1)\setminus [(T_1\times
[0,1))\cup(T_2\times [0,1))]$, and $\Sigma'_0$ a properly embedded
convex surface in $\Sigma'\times S^1$ isotopic to $\Sigma'\times
\{\text{pt}\}$ that has Legendrian boundary intersecting the
dividing set of $\partial \Sigma'\times S^1$ efficiently. Since
$\xi|_{\Sigma'\times S^1}$ is appropriate, the dividing set of
$\Sigma'_0$ consists of three properly embedded arc, each of which
connects a different pair of boundary components of $\Sigma'_0$.
Up to isotopy relative to $\partial\Sigma'_0$, there are
infinitely many such multi-arcs on $\Sigma'_0$. But, up to isotopy
of $\Sigma'_0$ which leaves $\partial\Sigma'_0$ invariant, there
are only two, each represented by a diagram in Figure 1 below.
Such an isotopy of $\Sigma'_0$ extend to an isotopy of
$\Sigma'\times S^1$ which, when restricted to a component of
$\partial\Sigma'\times S^1$, is a horizontal rotation. Thus, up to
isotopy of $\Sigma'\times S^1$, which, when restricted to any of
the components of $\partial\Sigma'\times S^1$, is a horizontal
rotation, there are only two appropriate contact structures on
$\Sigma'\times S^1$. Now let $\Phi_t$ be such an isotopy of
$\Sigma'\times S^1$ changing $\xi|_{\Sigma'\times S^1}$ to one of
the two standard appropriate contact structures. We extend
$\Phi_t$ to an isotopy $\widetilde{\Phi}_t$ of $\Sigma\times S^1$,
which fixes a neighborhood of $T_1\cup T_2$, and leaves $T_1\times
I$, $T_2\times I$ and $\Sigma'\times S^1$ invariant. Note that the
relative Euler class of $\xi|_{T_i\times I}$ is $(2k-n, 0)^T$,
where $k$ is the number of positive basic slices contained in
$(T^2\times I,\eta)$, and is invariant under
$\widetilde{\Phi}_t|_{T_i\times I}$. So $\xi|_{T_i\times I}$ and
$\widetilde{\Phi}_{1\ast}(\xi)|_{T_i\times I}$ have the same
relative Euler class, and are both continued fraction blocks
satisfying the same boundary condition. According to the
classification of tight contact structures on $T^2\times I$,
$\xi|_{T_i\times I}$ and
$\widetilde{\Phi}_{1\ast}(\xi)|_{T_i\times I}$ are isotopic
relative to boundary. So $\widetilde{\Phi}_{1\ast}(\xi)$ satisfies
the conditions given in the lemma, and is of one of the two
standard form. Thus, up to isotopy fixing $T_1$, $T_2$ and leaving
$T_3$ invariant, there are only two appropriate contact structures
on $\Sigma\times S^1$ satisfying the given conditions. Rotating
the diagram on the left of Figure 1 by $180^{\circ}$ induces a
self-diffeomorphism of $\Sigma\times S^1$ mapping $T_1$ to $T_2$
and changing the dividing set of $\Sigma'_0$ on the left of Figure
1 to the one on the right. So this self-diffeomorphism is isotopic
to a contactomorphism between the two standard appropriate contact
structures on $\Sigma\times S^1$. Hence, up to contactomorphism,
there is only one such appropriate contact structure on
$\Sigma\times S^1$. Thus, we can show the existence of an annulus
with the required properties by explicitly constructing such an
annulus in a model contact structure on $\Sigma\times S^1$ which
satisfies the given conditions.

\begin{figure}[h]
  \includegraphics[width=3in]{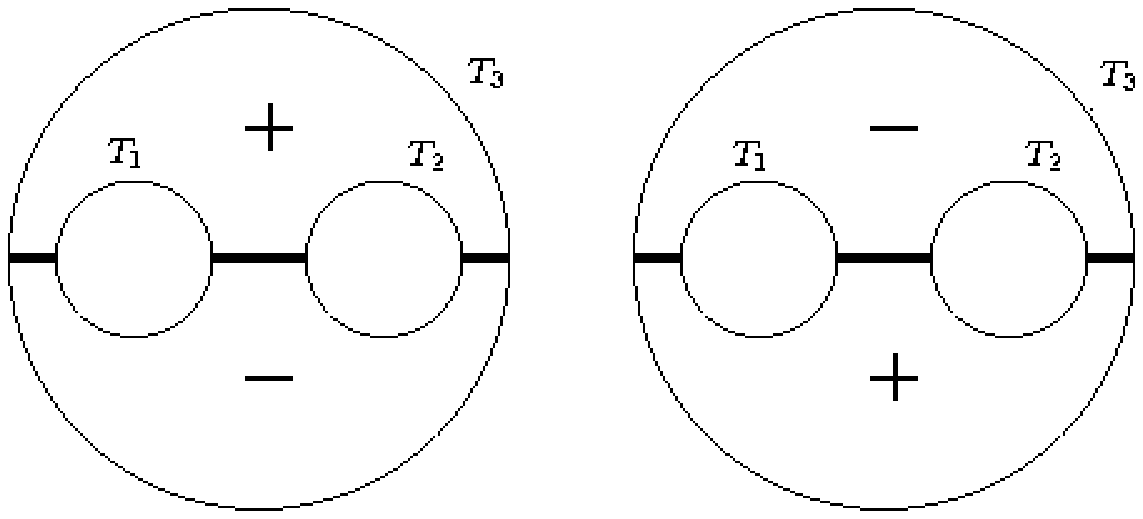}
  \caption{}\label{pants dividing}
\end{figure}

Consider the minimal twisting tight contact structure $\eta$ on
the thickened torus $T^2\times I$. Note that the vertical
Legendrian rulings of $T^2\times\{0\}$ intersect its dividing
curves efficiently. Without loss of generality, we assume that
$T^2\times \{1\}$ has horizontal Legendrian rulings and two
vertical Legendrian dividings. We further assume that, for a small
$\varepsilon>0$, $\eta|_{T^2\times[0,\varepsilon]}$ is invariant
in the $I$ direction. This is legitimate since $T^2\times \{0\}$
is convex. So $T^2\times \{\frac{\varepsilon}{2}\}$ is also a
convex torus with vertical Legendrian rulings and dividing curves
of slope $-\frac{1}{n}$. Let $L$ be a Legendrian ruling of
$T^2\times \{\frac{\varepsilon}{2}\}$. Since the twisting number
of $\eta|_L$ relative to the framing given by $T^2\times
\{\frac{\varepsilon}{2}\}$ is $-n$, we can find a standard
neighborhood $U$ of $L$ in $T^2\times(0,\varepsilon)$ such that
$\partial U$ is convex with vertical Legendrian ruling and two
dividing curves of slope $-\frac{1}{n}$. Now, we set $\Sigma\times
S^1= (T^2\times I)\setminus U$, where $T_1=T^2\times\{0\}$,
$T_2=\partial U$ and $T_3=-T^2\times\{1\}$, and let
$\xi=\eta|_{\Sigma\times S^1}$. Since $\eta$ is tight, so is
$\xi$. And there are no embedded thickened tori in $\Sigma\times
S^1$ with convex boundary isotopic to $T_2$ and $I$-twisting at
least $\pi$. Otherwise, $L$ would have an overtwisted neighborhood
in $T^2\times I$, which contradicts the tightness of $\eta$. Also,
since the $I$-twisting of $\eta$ is less than $\pi$, there exists
no embedded thickened tori in $\Sigma\times S^1$ with convex
boundary isotopic to $T_1$ or $T_3$ and $I$-twisting at least
$\pi$. Thus, $\xi$ is appropriate. Now we choose a vertical convex
annulus $A_1$ in $\Sigma\times S^1$ connecting a Legendrian ruling
of $T_1$ to a Legendrian divide of $T_3$, and a vertical convex
annulus $A_2$ in $\Sigma\times S^1$ connecting a Legendrian ruling
of $T_2$ to the other Legendrian divide of $T_3$ such that
$(T_1\cup A_1)\cap(T_2\cup A_2)=\phi$. The dividing set of $A_i$
consists of $n$ arcs starting and ending on $A_i\cap T_i$. We can
find a collar neighborhood $T_i\times I$ of $T_i$, for which
$T_i\times\{0\}=T_i$ and $T_i\times\{1\}$ is convex with dividing
set consisting of two circles of slope $\infty$, by isotoping
$T_i$ to engulf all the dividing curves of $A_i$ through bypass
adding. Since bypass adding can be done in a small neighborhood of
the original surface and the bypass, we can make $T_1\times I$ and
$T_2\times I$ mutually disjoint and disjoint from $T_3$. Note that
both $T_1\times I$ and $T_2\times I$ are minimal twisting. So they
are continued fraction blocks satisfying the boundary conditions
specified above. Let $k_i$ be the number of positive slices in
$T_i\times I$, and $B_i=A_i\cap (T_i\times I)$. Then
$2k_i-n=\chi((B_i)_+)-\chi((B_i)_-)=\chi((A_i)_+)-\chi((A_i)_-)$.
But $\chi((A_1)_+)-\chi((A_1)_-)=2k-n$, where $k$ is the number of
positive basic slices contained in $(T^2\times I,\eta)$. So
$k_1=k$. And, since $\eta|_{T^2\times(0,\varepsilon)}$ is
$I$-invariant, we can extend $A_2$ to a vertical annulus
$\widetilde{A}_2$ in $T^2\times I$ starting at a Legendrian ruling
of $T_1$ and such that
$\chi((\widetilde{A}_2)_+)-\chi((\widetilde{A}_2)_-)=\chi((A_2)_+)-\chi((A_2)_-)$.
Clearly,
$2k-n=\chi((\widetilde{A}_2)_+)-\chi((\widetilde{A}_2)_-)$. So
$k_2=k$. Thus, $k_1=k_2=k$. But the isotopy type of a continued
fraction block is determined by the number of positive slices in
it. Thus, $\xi|_{T_1\times I}$, $\xi|_{T_2\times I}$ and $\eta$
are isotopic relative to boundary. So our $(\Sigma\times S^1,
\xi)$ is indeed a legitimate model. Now we connect a Legendrian
ruling of $T_1$ and a Legendrian ruling of $T_2$ by a vertical
convex annulus $A$ which is contained in $(T^2\times
[0,\varepsilon])\setminus U$. Then $\partial A$ intersects the
dividing sets of $T_1$ and $T_2$ efficiently. If $A$ has
$\partial$-parallel diving curves, then $(T^2\times
[0,\varepsilon])$ has non-zero $I$-twisting, which contradicts our
choice of the slice $(T^2\times [0,\varepsilon])$. Thus, $A$ has
no $\partial$-parallel diving curves.
\end{proof}

Now we are in the position to prove Theorem \ref{-2}.

\vspace{.4cm}

\noindent\textit{Proof of Theorem \ref{-2}.} Let
$M=M(\frac{q_1}{p_1},\frac{q_2}{p_2},\frac{q_3}{p_3})$ be a small
Seifert space with $e_0(M)\leq-2$. Without loss of generality, we
assume that $p_1,p_2,p_3>1$, $0<q_1<p_1$, and $q_2,q_3<0$. For
$i=1,2,3$, let $V_i=D^2\times S^1$, and identify $\partial V_i$
with $\bf{R}^2/\bf{Z}^2$ by identifying a meridian $\partial D^2
\times \{\text{pt}\}$ with $(1,0)^T$ and a longitude
$\{\text{pt}\}\times S^1$ with $(0,1)^T$. Choose $u_i, v_i \in
\bf{Z}$ such that $p_iv_i+q_iu_i=1$ for $i=1,2,3$. Define an
orientation preserving diffeomorphism $\varphi_i:\partial
V_i\rightarrow T_i$ by
\[
\varphi_i ~=~
\left(%
\begin{array}{cc}
  p_i & u_i \\
  -q_i & v_i \\
\end{array}%
\right),
\]
where $T_i$ and the coordinates on it are defined above. Then
\[
M=M(\frac{q_1}{p_1},\frac{q_2}{p_2},\frac{q_3}{p_3})\cong(\Sigma
 \times S^1)\cup_{(\varphi_1\cup\varphi_2\cup\varphi_3)}(V_1\cup V_2\cup
 V_3).
\]

Assume that $\xi$ is a tight contact structure on $M$ for which
there exits a Legendrian vertical circle $L$ in $M$ with twisting
number $t(L)=0$. We first isotope $\xi$ to make
$L=\{\text{pt}\}\times S^1\subset \Sigma \times S^1$, and each
$V_i$ a standard neighborhood of a Legendrian circle $L_i$
isotopic to the $\frac{q_i}{p_i}$-singular fiber with twisting
number $n_i<0$, i.e., $\partial V_i$ is convex with two dividing
curves each of which has slope $\frac{1}{n_i}$ when measured in
the coordinates of $\partial V_i$ given above. Let $s_i$ be the
slope of the dividing curves of $T_i=\partial V_i$ measured in the
coordinates of $T_i$. Then we have that
\[
s_i ~=~ \frac{-n_iq_i+v_i}{n_ip_i+u_i} ~=~
-\frac{q_i}{p_i}+\frac{1}{p_i(n_ip_i+u_i)}.
\]
From our choice of $p_i$ and $q_i$, one can see that $-1 \leq s_1
\leq 0$ and $0\leq s_2,s_3 < \infty$. Now, without affecting the
properties of $L$ and $V_i$ asserted above, we can further isotope
the contact structure $\xi$ to make the Legendrian rulings of
$T_i$ to have slope $\infty$ when measured in the coordinates of
$T_i$.

Pick a Legendrian ruling $\widetilde{L}_i$ on each $T_i$, and
connect $L$ to $\widetilde{L}_i$ by a vertical convex annulus
$A_i$ such that $A_i\cap A_j=L$ when $i\neq j$. Let $\Gamma_{A_i}$
be the dividing set of $A_i$. Since $A_i$ gives the canonical
framing $Fr$ of $L$, we know that the twisting number of $\xi|_L$
relative to $TA_i|_L$ is $0$. This means that $\Gamma_{A_i}\cap
L=\phi$. But $\Gamma_{A_i}\cap \widetilde{L}_i\neq\phi$. There are
dividing curves of $A_i$ starting and ending on $\widetilde{L}_i$.
According to Lemma 3.15 of \cite{H1}, we can find an embedded
minimal twisting slice $T_i\times I$ in $\Sigma \times S^1$, for
which $T_i\times\{0\}=T_i$, $T_i\times\{1\}$ is convex with two
vertical dividing curves, by isotoping $T_i$ to engulf all the
dividing curves of $A_i$ starting and ending on $\widetilde{L}_i$
through bypass adding. Since bypass adding can be done in a small
neighborhood of the bypass and the original surface, and the
bypasses from different $A_i$'s are mutually disjoint, we can make
$T_i\times I$'s pairwise disjoint. By Corollary 4.16 of \cite{H1},
we can find a convex torus in $T_i\times(0,1)$ isotopic to $T_i$
that has two dividing curves of the slope $-1$. Without loss of
generality, we assume that this torus is
$T_i\times\{\frac{1}{2}\}$. Moreover, for $i=2,3$, we can find
another convex torus, say $T_i\times\{\frac{1}{4}\}$, in
$T_i\times(0,\frac{1}{2})$ isotopic to $T_i$ with two dividing
curves of slope $0$.

\begin{figure}[h]
  \includegraphics[width=3in]{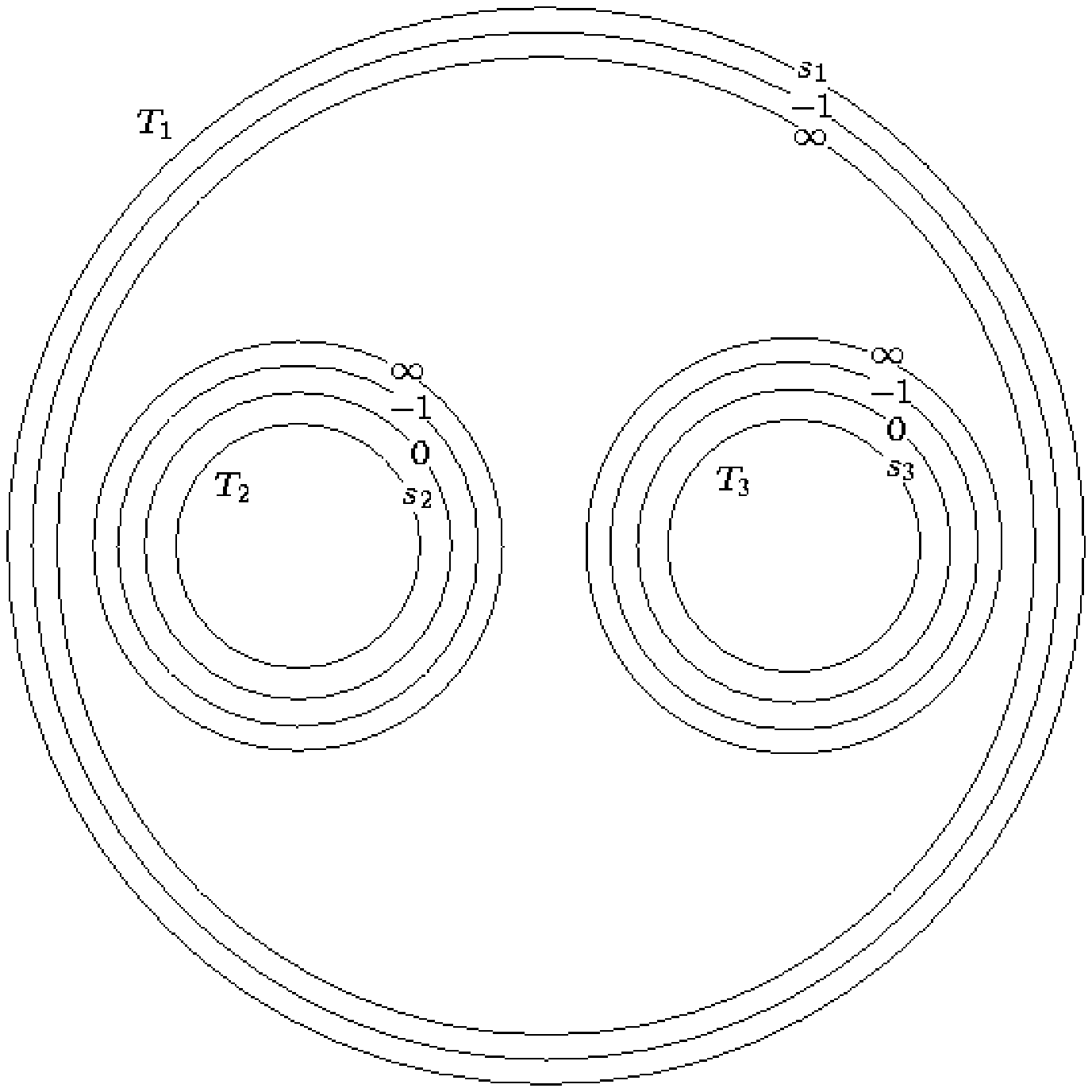}
  \caption{}\label{1}
\end{figure}

Since the slice $T_i\times I$ is minimal twisting, so is any of
its sub-slices. Let's consider the thickened tori
$T_i\times[\frac{1}{2},1]$. All of these have the same boundary
condition, and are minimal twisting. There are only two such tight
contact structures up to isotopy relative to boundary. So two of
these have to be isotopic relative to boundary. There are $3$
cases.

\textit{Case 1.} $T_1\times[\frac{1}{2},1]$ and
$T_2\times[\frac{1}{2},1]$ are isotopic. We apply Lemma \ref{GS36}
to
\[
\Sigma' \times S^1=(\Sigma \times S^1)\setminus
(T_1\times[0,\frac{1}{2}) \cup T_2\times[0,\frac{1}{2}) \cup
T_3\times[0,1)).
\]
Then there exists a vertical convex annulus $A$ connecting
$T_1\times\{\frac{1}{2}\}$ and $T_2\times\{\frac{1}{2}\}$ with no
$\partial$-parallel dividing curves that has Legendrian boundary
intersecting the dividing sets of these tori efficiently. We can
extend $A$ across $T_2\times[\frac{1}{4},\frac{1}{2}]$ to a convex
annulus $\widetilde{A}$ connecting $T_1\times\{\frac{1}{2}\}$ and
$T_2\times\{\frac{1}{4}\}$ with Legendrian boundary intersecting
the dividing sets of these two tori efficiently. Since
$T_2\times[\frac{1}{4},\frac{1}{2}]$ is minimal twisting,
$\widetilde{A}\setminus A$ has no $\partial$-parallel dividing
curves. Thus, $\widetilde{A}$ has no $\partial$-parallel dividing
curves either. Cut $(\Sigma \times S^1)\setminus
(T_1\times[0,\frac{1}{2}) \cup T_2\times[0,\frac{1}{4}) \cup
T_3\times[0,1))$ along $\widetilde{A}$, and round the edges. We
get a thickened torus $T_3\times[1,2]$ embedded in $\Sigma \times
S^1$ with convex boundary. The dividing set of $T_3\times\{2\}$
consists two circles of slope $0$. Now we can see that the
thickened torus $T_3\times[0,2]$ has $I$-twisting at least $\pi$
since the dividing curves of $T_3\times\{\frac{1}{4}\}$ and
$T_3\times\{2\}$ have slope $0$ and those of $T_3\times\{1\}$ have
slope $\infty$. Thus, $\Sigma \times S^1$ is inappropriate. This
is a contradiction.

\textit{Case 2.} $T_1\times[\frac{1}{2},1]$ and
$T_3\times[\frac{1}{2},1]$ are isotopic. The proof for this case
is identical to that of Case 1 except for interchanging the
subindexes $2$ and $3$.

\textit{Case 3.} $T_2\times[\frac{1}{2},1]$ and
$T_3\times[\frac{1}{2},1]$ are isotopic. Similar to Case 1, we can
find a vertical convex annulus $B$ connecting
$T_2\times\{\frac{1}{2}\}$ and $T_3\times\{\frac{1}{2}\}$ with no
$\partial$-parallel dividing curves that has Legendrian boundary
intersecting the dividing sets of these tori efficiently. Extend
$B$ across $T_2\times[\frac{1}{4},\frac{1}{2}]$ and
$T_3\times[\frac{1}{4},\frac{1}{2}]$ to a convex annulus
$\widetilde{B}$ connecting $T_2\times\{\frac{1}{4}\}$ and
$T_3\times\{\frac{1}{4}\}$ with Legendrian boundary intersecting
the dividing sets of these two tori efficiently. For reasons
similar to above, neither component of $\widetilde{B}\setminus B$
has $\partial$-parallel dividing curves. Thus, $\widetilde{B}$ has
no $\partial$-parallel dividing curves. Cut $(\Sigma \times
S^1)\setminus (T_1\times[0,1) \cup T_2\times[0,\frac{1}{4}) \cup
T_3\times[0,\frac{1}{4}))$ along $\widetilde{B}$, and round the
edges. We get a thickened torus $T_1\times[1,2]$ embedded in
$\Sigma \times S^1$ with convex boundary. The dividing set of
$T_1\times\{2\}$ consists two circles of slope $-1$. Now we can
see that the thickened torus $T_1\times[0,2]$ has $I$-twisting at
least $\pi$ since the dividing curves of
$T_1\times\{\frac{1}{2}\}$ and $T_1\times\{2\}$ have slope $-1$
and those of $T_3\times\{1\}$ have slope $\infty$. Thus, $\Sigma
\times S^1$ is inappropriate. This is again a contradiction.

Thus, $M=M(\frac{q_1}{p_1},\frac{q_2}{p_2},\frac{q_3}{p_3})$
admits no tight contact structures for which there exists a
Legendrian vertical circle with twisting number $0$. \hfill $\Box$

\vspace{.4cm}

\section{The $e_0=-1$ Case}

Since part (1) of Theorem \ref{-1} is already proved, we will
concentrate on parts (2) and (3) of Theorem \ref{-1}. We will
refine the method used in the $e_0\leq-2$ case to prove these
results. Lemmata \ref {sign}, \ref{GS37} below and Lemma
\ref{GS36} from last section will be the main technical tools used
in the proof.

\begin{lemma}\label{sign}
Let $\xi$ be an appropriate contact structure on $\Sigma\times
S^1$. Suppose that $-\partial\Sigma\times S^1=T_1+T_2+T_3$ is
convex and such that each of $T_1$ and $T_2$ has two dividing
curves of slope $-1$, and $T_3$ has two horizontal dividing
curves. Assume that there are pairwise disjoint collar
neighborhoods $T_i\times I$ of $T_i$ in $\Sigma\times S^1$ for
$i=1,2,3$, such that $T_i\times\{0\}=T_i$ and $T_i\times\{1\}$ is
convex with two vertical dividing curves. Then $(T_1\times
I,\xi|_{T_1\times I})$, $(T_2\times I,\xi|_{T_2\times I})$ and
$(T_3\times I,\xi|_{T_3\times I})$ are all basic slices, and the
signs of these basic slices can not be all the same, where the
sign convention of $(T_i\times I,\xi|_{T_i\times I})$ is given by
choosing the vector associated with $T_i\times \{1\}$ to be
$(0,1)^T$.
\end{lemma}

\begin{proof}
Since $\xi$ is appropriate, each $(T_i\times I,\xi|_{T_i\times
I})$ is minimal twisting. From the boundary condition of these
slices, we can see these are all basic slices. Assume that all
these basic slices have the same sign. Then we have that
$(T_1\times I,\xi|_{T_1\times I})$ and $(T_2\times
I,\xi|_{T_2\times I})$ are isotopic relative to boundary. We
isotope $T_1$ and $T_2$ slightly so that they have vertical
Legendrian rulings. By Lemma \ref{GS36}, we can then find a
properly embedded convex vertical annulus $A$ disjoint from $T_3
\times I$ with no $\partial$-parallel dividing curves, whose
Legendrian boundary $\partial A=(A\cap T_1)\cup (A\cap T_2)$
intersects the dividing sets of $T_1$ and $T_2$ efficiently. Cut
$\Sigma\times S^1$ open along $A$, we get a thickened torus
$T_3\times [0,2]$ such that each of $T_3\times\{0\}$,
$T_3\times\{1\}$ and $T_3\times\{2\}$ is convex with two dividing
curves, and the slopes of the dividing curves are $0$, $\infty$
and $1$, respectively. Note that the slice $(T_3\times
[1,2],\xi|_{T_3\times [1,2]})$ has the sign different from that of
$(T_1\times I,\xi|_{T_1\times I})$, and the slice $(T_3\times
[0,1],\xi|_{T_3\times [0,1]})$ has the same sign as that of
$(T_1\times I,\xi|_{T_1\times I})$. So $\xi_{T_3\times [0,2]}$ is
overtwisted. This is a contradiction. Thus, the signs of the basic
slices $(T_1\times I,\xi|_{T_1\times I})$, $(T_2\times
I,\xi|_{T_2\times I})$ and $(T_3\times I,\xi|_{T_3\times I})$ can
not be all the same.
\end{proof}

The following lemma is a special case of Lemma 37 of \cite{GS}.
Its proof is quite similar to that of Lemma \ref{GS36} (\cite{GS},
Lemma 36). We will only give a sketch of it.

\begin{lemma}[\cite{GS}, Lemma 37]\label{GS37} Let $\xi$
be an appropriate contact structure on $\Sigma\times S^1$. Suppose
that $-\partial\Sigma\times S^1=T_1+T_2+T_3$ is convex and such
that $T_1$ has vertical Legendrian rulings and two dividing curves
of slope $-\frac{1}{n}$, where $n\in{\bf Z}^{>0}$, $T_2$ has
vertical Legendrian rulings and two dividing curves of slope
$\frac{1}{n}$, and $T_3$ has two vertical dividing curves. Let
$T_1\times I$ and $T_2\times I$ be collar neighborhoods of $T_1$
and $T_2$ that are mutually disjoint and disjoint from $T_3$, and
such that, for $i=1,2$, $T_i\times\{0\}=T_i$ and $T_i\times\{1\}$
is convex with dividing set consisting of two  vertical circles.
If basic slices of $(T_1\times I,\xi|_{T_1\times I})$ and
$(T_2\times I,\xi|_{T_2\times I})$ are all of the same sign, then
there exists a properly embedded convex vertical annulus $A$ with
no $\partial$-parallel dividing curves, whose Legendrian boundary
$\partial A=(A\cap T_1)\cup (A\cap T_2)$ intersects the dividing
sets of $T_1$ and $T_2$ efficiently.
\end{lemma}

\noindent\textit{Sketch of proof.} Similar to the proof of Lemma
\ref{GS36}, we can show that, if we prescribe the sign of the
basic slices of $(T_1\times I,\xi|_{T_1\times I})$ and $(T_2\times
I,\xi|_{T_2\times I})$, then, up to isotopy that fixes $T_1$,
$T_2$ and leaves $T_3$ invariant, there are at most two
appropriate contact structures on $\Sigma\times S^1$ that satisfy
the given conditions, each of which corresponds to one of the two
diagrams in Figure \ref{pants dividing}. Since the two layers
$T_1\times I$ and $T_2\times I$ are not contactomorphic, we can
not find a contactomorphism between these two possible appropriate
contact structures as before. Instead, we will construct an
appropriate contact structure corresponding to each of these two
diagrams, and show that each of these admits an annulus with the
required properties.

\begin{figure}[h]
  \includegraphics[width=2in]{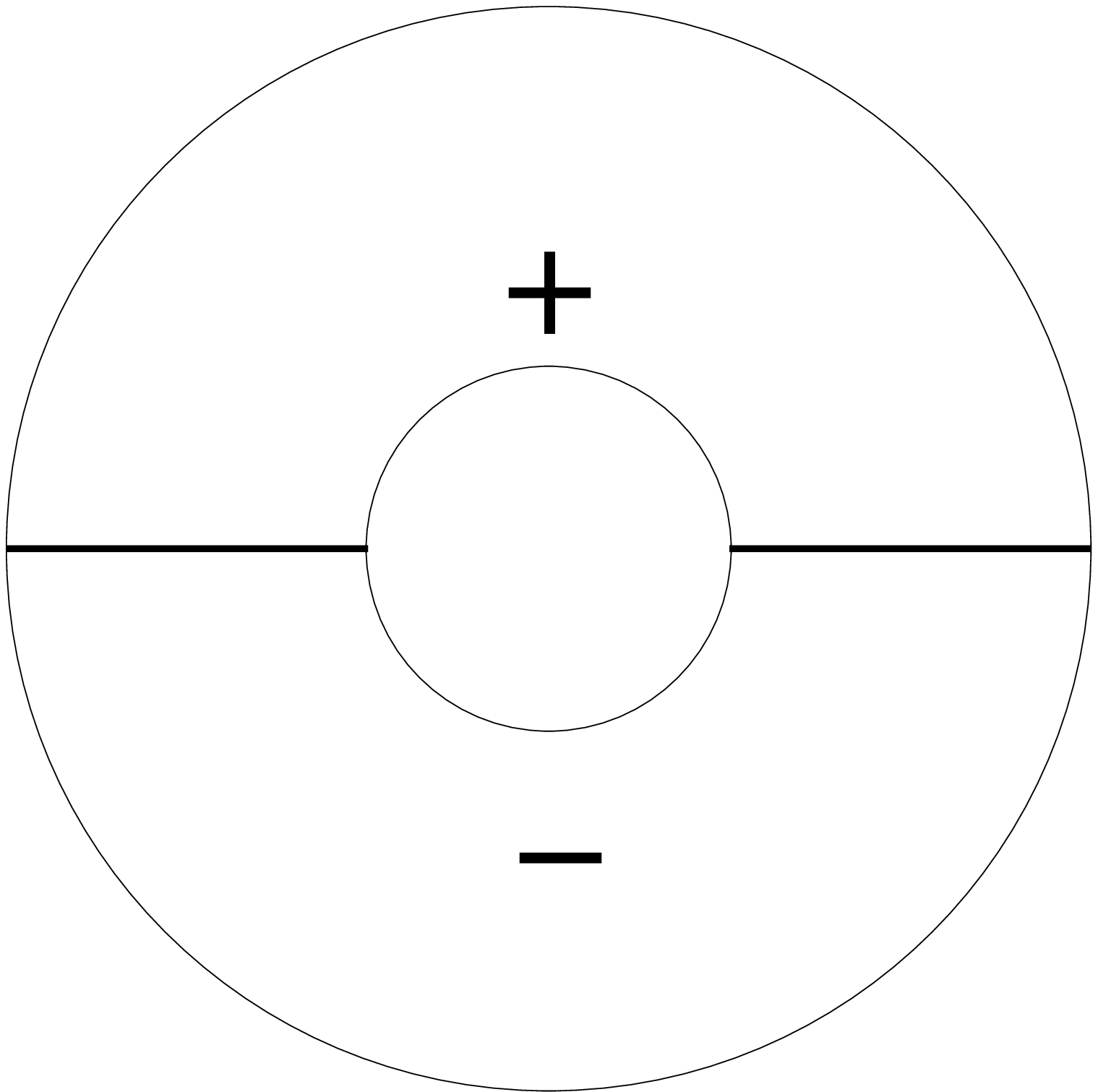}
  \caption{}\label{band dividing}
\end{figure}

Now consider the tight contact thickened torus $(T_2\times
I,\xi|_{T_2\times I})$. Like in the proof of Lemma \ref{GS36}, we
can construct an appropriate contact structure on $\Sigma\times
S^1$ satisfying the conditions in the lemma that admits an annulus
$A$ with the required properties by "digging out" a vertical
Legendrian ruling of a torus in an $I$-invariant neighborhood of
$T_2 \times \{0\}$ parallel to the the boundary. Indeed, both of
the possible appropriate contact structures can be constructed
this way. To see that, we isotope $T_2 \times \{0\}$ and
$T_2\times \{1\}$ lightly to $T'_2$ and $T'_3$ with the same
dividing curves and horizontal Legendrian rulings. Then connect a
Legendrian ruling of $T'_2$ and a Legendrian ruling of $T'_3$ by a
horizontal convex annulus $B$. The dividing curves of $B$ is given
in Figure \ref{band dividing}. We can choose the vertical
Legendrian ruling to be dug out to intersect one of the two
dividing curves of $B$. These two choices of the vertical
Legendrian ruling correspond to the two possible layout of the
dividing curves on $\Sigma'_0$ in Figure \ref{pants dividing},
and, hence, gives the two possible appropriate contact structures
on $\Sigma\times S^1$ satisfying the given conditions. \hfill
$\Box$

\vspace{.4cm}

\noindent\textit{Proof of (2) and (3) of Theorem \ref{-1}.} Let
$M=M(\frac{q_1}{p_1},\frac{q_2}{p_2},\frac{q_3}{p_3})$ be a small
Seifert space such that $0<q_1<p_1$, $0<q_2<p_2$ and $-p_3<q_3<0$.
For $i=1,2,3$, let $V_i=D^2\times S^1$, and identify $\partial
V_i$ with $\bf{R}^2/\bf{Z}^2$ by identifying a meridian $\partial
D^2 \times \{\text{pt}\}$ with $(1,0)^T$ and a longitude
$\{\text{pt}\}\times S^1$ with $(0,1)^T$. Choose $u_i, v_i \in
\bf{Z}$ such that $0<u_i<p_i$ and $p_iv_i+q_iu_i=1$ for $i=1,2,3$.
Define an orientation preserving diffeomorphism
$\varphi_i:\partial V_i\rightarrow T_i$ by
\[
\varphi_i ~=~
\left(%
\begin{array}{cc}
  p_i & u_i \\
  -q_i & v_i \\
\end{array}%
\right),
\]
where $T_i$ and the coordinates on it are defined above. Then
\[
M=M(\frac{q_1}{p_1},\frac{q_2}{p_2},\frac{q_3}{p_3})\cong(\Sigma
 \times S^1)\cup_{(\varphi_1\cup\varphi_2\cup\varphi_3)}(V_1\cup V_2\cup
 V_3).
\]

Assume that $\xi$ is a tight contact structure on $M$ for which
there exits a Legendrian vertical circle $L$ in $M$ with twisting
number $t(L)=0$. We isotope $\xi$ to make $L=\{\text{pt}\}\times
S^1\subset \Sigma \times S^1$, and each $V_i$ a standard
neighborhood of a Legendrian circle $L_i$ isotopic to the $i$-th
singular fiber with twisting number $n_i<0$, i.e., $\partial V_i$
is convex with two dividing curves each of which has slope
$\frac{1}{n_i}$ when measured in the coordinates of $\partial V_i$
given above. Let $s_i$ be the slope of the dividing curves of
$T_i=\partial V_i$ measured in the coordinates of $T_i$. Then we
have that
\[
s_i ~=~ \frac{-n_iq_i+v_i}{n_ip_i+u_i} ~=~
-\frac{q_i}{p_i}+\frac{1}{p_i(n_ip_i+u_i)}.
\]
Then $-1 \leq s_1,s_2 \leq 0$ and $0\leq s_3 < 1$. As before, we
can find pairwise disjoint collar neighborhoods  $T_i\times I$'s
in $\Sigma \times S^1$ of $T_i$'s, such that $T_i\times\{0\}=T_i$
and $T_i\times\{1\}$ is convex with dividing set consisting of two
vertical circles.

We now prove part (2).

Assume that $q_3=-1$, $\frac{q_1}{p_1}<\frac{1}{2p_3-1}$ and
$\frac{q_2}{p_2}<\frac{1}{2p_3}$. By choosing $n_i\ll-1$, we can
make $-\frac{1}{2p_3-1}<s_1<-\frac{q_1}{p_1}$,
$-\frac{1}{2p_3}<s_2<-\frac{q_2}{p_2}$ and
$\frac{1}{p_3+1}<s_3<\frac{1}{p_3}$. So there is a convex torus in
$T_i \times I$ parallel to the boundary, say $T'_i=T_i \times
\{\frac{1}{2}\}$, that has two dividing curves of slope
$-\frac{1}{2p_3-1}$, $-\frac{1}{2p_3}$ and $\frac{1}{p_3+1}$ for
$i=1,2 ~\text{and}~3$, respectively. Let's consider the layers
$T_i \times [\frac{1}{2},1]$. $T_1 \times [\frac{1}{2},1]$ is a
continued fraction block consisting of $2p_3-1$ basic slices. $T_2
\times [\frac{1}{2},1]$ is a continued fraction block consisting
of $2p_3$ basic slices. $T_3 \times [\frac{1}{2},1]$ consists of
$2$ continued fraction blocks, each of which is a basic slice. We
can find a convex torus $T''_i=T_i \times \{\frac{3}{4}\}$ in $T_i
\times [\frac{1}{2},1]$ parallel to bounary with two dividing
curves of slope $-1$ for $i=1,2$, and $0$ for $i=3$.

Let the sign of the basic slice $T_3\times[\frac{3}{4},1]$ be
$\sigma\in\{+,-\}$. Note that, when $q_3=-1$, then diffeomorphism
$\varphi_3:\partial V_3\rightarrow T_3$ is given by
\[
\varphi_3 ~=~
\left(%
\begin{array}{cc}
  p_3 & p_3-1 \\
  1 & 1 \\
\end{array}%
\right).
\]
So the slopes $0$ and $\frac{1}{p_3+1}$ of the dividing sets of
$T''_3$ and $T'_3$ correspond to twisting numbers $-1$ and $-2$ of
Legendrian circles isotopic to the $-\frac{1}{p_3}$-singular
fiber. And the basic slice $T_3 \times [\frac{1}{2},\frac{3}{4}]$
corresponds to a stabilization of a Legendrian circle isotopic to
the $-\frac{1}{p_3}$-singular fiber. Since we can freely choose
the sign of such a stabilization, we can make the sign of the
basic slice $T_3 \times [\frac{1}{2},\frac{3}{4}]$ to be $\sigma$,
too.

According to Lemma \ref{sign}, the sign of the basic slices
$T_i\times[\frac{3}{4},1]$ can not be all the same. Note that we
can shuffle the signs of basic slices in a continued fraction
block. So at least one of $T_1 \times [\frac{1}{2},1]$ and $T_2
\times [\frac{1}{2},1]$ consists only of basic slices of sign
$-\sigma$.

\textit{Case 1.} Assume that all the basic slices in $T_1 \times
[\frac{1}{2},1]$ are of the sign $-\sigma$. If $T_2 \times
[\frac{1}{2},1]$ contains $p_3$ basic slices of the sign
$-\sigma$, then we shuffle these signs to the $p_3$ slices closest
to $T_2\times\{1\}$. Now consider the thickened tori $T_1 \times
[\frac{5}{8},1]$ and $T_2 \times [\frac{5}{8},1]$ formed by the
unions the $p_3$ basic slices closest to $T_1\times\{1\}$ and
$T_2\times\{1\}$ in $T_1 \times I$ and $T_2 \times I$,
respectively. Remove from $M$ the open solid tori bounded by $T_1
\times\{\frac{5}{8}\}$, $T_2 \times\{\frac{5}{8}\}$ and $T_3
\times\{1\}$. We apply Lemma \ref{GS36} to the resulting $\Sigma
\times S^1$ and the thickened tori $T_1 \times [\frac{5}{8},1]$
and $T_2 \times [\frac{5}{8},1]$. Then there exists a properly
embedded convex vertical annulus $A$ in $\Sigma \times S^1$ with
no $\partial$-parallel dividing curves, whose Legendrian boundary
$\partial A=(A\cap (T_1 \times\{\frac{5}{8}\}))\cup (A\cap (T_2
\times\{\frac{5}{8}\}))$ intersects the dividing sets of $T_1
\times\{\frac{5}{8}\}$ and $T_2 \times\{\frac{5}{8}\}$
efficiently. Cutting $\Sigma \times S^1$ open along $A$ and round
the edges, we get a torus convex $\widetilde{T}_3$ isotopic to
$T_3$ with two dividing curves of slope $\frac{1}{p_3}$. This
means there exists a thickening $\widetilde{V}_3$ of $V_3$ with
convex boundary $\partial\widetilde{V}_3$ that has two dividing
curves isotopic to a meridian. Then
$\xi|_{\partial\widetilde{V}_3}$ is overtwisted. This contradicts
the tightness of $\xi$.

If $T_2 \times [\frac{1}{2},1]$ contains $p_3+1$ basic slices of
the sign $\sigma$, then we shuffle all these signs to the $p_3+1$
slices closest to $T_2 \times \{1\}$. Let $T_2 \times
[\frac{5}{8},1]$ be the union of these $p_3+1$ basic slices.
Remove from $M$ the open solid tori bounded by $T_1\times\{1\}$,
$T_2 \times \{\frac{5}{8}\}$ and $T_3 \times \{\frac{1}{2}\}$.
Apply Lemma \ref{GS37} to the resulting $\Sigma \times S^1$ and
the thickened tori $T_2 \times [\frac{5}{8},1]$ and $T_3 \times
[\frac{1}{2},1]$. Then there exists a properly embedded convex
vertical annulus $A$ in $\Sigma \times S^1$ with no
$\partial$-parallel dividing curves, whose Legendrian boundary
$\partial A=(A\cap (T_2 \times\{\frac{5}{8}\}))\cup (A\cap (T_3
\times\{\frac{1}{2}\}))$ intersects the dividing sets of $T_2
\times\{\frac{5}{8}\}$ and $T_3 \times\{\frac{1}{2}\}$
efficiently. Cutting $\Sigma \times S^1$ open along $A$, we get a
thickened torus $T_1\times[1,2]$ embedded in $\Sigma \times S^1$
that has convex boundary such that $T_1\times\{1\}$ has vertical
dividing curves and $T_1\times\{2\}$ has dividing curves of slope
$-\frac{1}{p_3+1}$. Then the thickened torus
$T_1\times[\frac{1}{2},2]=(T_1\times[\frac{1}{2},1])\cup(T_1\times[1,2])$
has $I$-twisting at least $\pi$. This again contradicts the
tightness of $\xi$.

But $T_2 \times [\frac{1}{2},1]$ is a continued fraction block
consisting of $2p_3$ basic slices. So it contains either $p_3$
basic slices of the sign $-\sigma$ or $p_3+1$ basic slices of the
sign $\sigma$. So, the basic slices in $T_1 \times
[\frac{1}{2},1]$ can not be all of the sign $-\sigma$.

\textit{Case 2.} Assume that all the basic slices in $T_2 \times
[\frac{1}{2},1]$ are of the sign $-\sigma$. If $T_1 \times
[\frac{1}{2},1]$ contains either $p_3$ basic slices of the sign
$-\sigma$ or $p_3+1$ basic slices of the sign $\sigma$, then there
will be a contradiction just like in Case 1. So the only possible
scenario is that $T_1 \times [\frac{1}{2},1]$ contains $p_3-1$
basic slices of the sign $-\sigma$ and $p_3$ basic slices of the
sign $\sigma$. Now we shuffle all the $-\sigma$ signs in $T_1
\times [\frac{1}{2},1]$ to the $p_3-1$ basic slices closest to
$T_1 \times \{1\}$. Now let $T_1 \times [\frac{5}{8},1]$ and $T_2
\times [\frac{5}{8},1]$ be the unions the $p_3-1$ basic slices
closest to $T_1\times\{1\}$ and $T_2\times\{1\}$ in $T_1 \times I$
and $T_2 \times I$. Remove from $M$ the open solid tori bounded by
$T_1 \times \{\frac{5}{8}\}$, $T_2 \times \{\frac{5}{8}\}$ and
$T_3 \times \{1\}$, and apply Lemma \ref{GS36} to the resulting
$\Sigma \times S^1$ and the thickened tori $T_1 \times
[\frac{5}{8},1]$ and $T_2 \times [\frac{5}{8},1]$. Then there
exists a properly embedded convex vertical annulus $A$ in $\Sigma
\times S^1$ with no $\partial$-parallel dividing curves, whose
Legendrian boundary $\partial A=(A\cap (T_1
\times\{\frac{5}{8}\}))\cup (A\cap (T_2 \times\{\frac{5}{8}\}))$
intersects the dividing sets of $T_1 \times\{\frac{5}{8}\}$ and
$T_2 \times\{\frac{5}{8}\}$ efficiently. Cutting $\Sigma \times
S^1$ open along $A$ and round the edges, we get a torus convex
$\widetilde{T}_3$ isotopic to $T_3$ with two dividing curves of
slope $\frac{1}{p_3-1}$. This means we can thicken $V_3$ to a
standard neighborhood $\widetilde{V}_3$ of a Legendrian circle
isotopic to the $-\frac{1}{p_3}$-singular fiber with twisting
number $0$. Stabilize this Legendrian circle twice. We get a
thickened torus $\widetilde{T}_3\times[\frac{1}{2},2]$ such that
$\widetilde{T}_3\times\{2\}=\widetilde{T}_3$,
$\widetilde{T}_3\times\{\frac{3}{4}\}$ has two dividing curves of
slope $0$, and $\widetilde{T}_3\times\{\frac{1}{2}\}$ has two
dividing curves of slope $\frac{1}{p_3+1}$. Since we can choose
the sign of these stabilizations freely, we can make both basic
slices $\widetilde{T}_3\times[\frac{1}{2},\frac{3}{4}]$ and
$\widetilde{T}_3\times[\frac{3}{4},2]$ to have the sign $-\sigma$.
There exists a convex torus, say $\widetilde{T}_3\times\{1\}$, in
$\widetilde{T}_3\times[\frac{3}{4},2]$ parallel to boundary that
has two vertical dividing curves. Use
$\widetilde{T}_3\times\{1\}$, we can thicken $T_1 \times
[\frac{1}{2},\frac{5}{8}]$ to $\widetilde{T}_1 \times
[\frac{1}{2},1]$, such that $\widetilde{T}_1 \times
[\frac{1}{2},\frac{5}{8}]=T_1 \times [\frac{1}{2},\frac{5}{8}]$,
and $\widetilde{T}_1 \times \{1\}$ is convex with two vertical
dividing curves. Since the basic
$\widetilde{T}_3\times[\frac{3}{4},2]$ has the sign $-\sigma$, all
the basic slices in $\widetilde{T}_1 \times [\frac{5}{8},1]$ have
the sign $\sigma$. Also note that all the basic slices in
$\widetilde{T}_1 \times [\frac{1}{2},\frac{5}{8}]=T_1 \times
[\frac{1}{2},\frac{5}{8}]$ have the sign $\sigma$. So we are now
in a situation where the basic slices
$\widetilde{T}_3\times[\frac{1}{2},\frac{3}{4}]$ and
$\widetilde{T}_3\times[\frac{3}{4},1]$ both have the sign
$-\sigma$, and all the basic slices in $\widetilde{T}_1 \times
[\frac{1}{2},1]$ have the sign $\sigma$. After we thicken $T_2
\times [\frac{1}{2},\frac{5}{8}]$ to $\widetilde{T}_2 \times
[\frac{1}{2},1]$, where $\widetilde{T}_2 \times \{1\}$ is convex
with two vertical dividing curves, we are back to Case 1, which is
shown to be impossible. Thus, the basic slices in $T_2 \times
[\frac{1}{2},1]$ can not be all of the sign $-\sigma$ either.

But, as we mentioned above, one of $T_1 \times [\frac{1}{2},1]$
and $T_2 \times [\frac{1}{2},1]$ have to consist only of basic
slices of sign $-\sigma$. This is a contradiction. Thus, no such
$\xi$ exists on $M$, and, hence, we proved part (2) of Theorem
\ref{-1}.

It remains to prove part(3) now.

Assume that $q_1=q_2=1$ and $p_1,p_2>2m$, where
$m=-\lfloor\frac{p_3}{q_3}\rfloor$. By choosing $n_i\ll-1$, we can
make $-\frac{1}{2m}<s_1<-\frac{1}{p_1}$,
$-\frac{1}{2m}<s_2<-\frac{1}{p_2}$, and $0<s_3<-\frac{q_3}{p_3}$.
Similar to the proof of part (2), we can find convex a torus
$T'_i=T_i \times \{\frac{1}{2}\}$ in $T_i \times I$ parallel to
boundary with two dividing curves that have slope $-\frac{1}{2m}$
for $i=1,2$, and $0$ for $i=3$. Then each of $T_1 \times
[\frac{1}{2},1]$ and $T_2 \times [\frac{1}{2},1]$ is a continued
fraction block consists of $2m$ basic slices. And $T_3 \times
[\frac{1}{2},1]$ is a basic slice. Let the sign of the basic slice
$T_3\times[\frac{1}{2},1]$ be $\sigma\in\{+,-\}$. For reasons
similar to above, at least one of $T_1 \times [\frac{1}{2},1]$ and
$T_2 \times [\frac{1}{2},1]$ can not contain basic slices of the
sign $\sigma$. Without loss of generality, we assume that all
basic slices in $T_1 \times [\frac{1}{2},1]$ are of the sign
$-\sigma$. We now consider the signs of the basic slices in $T_2
\times [\frac{1}{2},1]$.

\textit{Case 1.} Assume that $T_2 \times [\frac{1}{2},1]$ contains
$m$ basic slices of the sign $-\sigma$. Then we shuffle these
signs to the $m$ basic slices in $T_2 \times [\frac{1}{2},1]$
closest to $T_2 \times \{1\}$. Denote by $T_1 \times
[\frac{3}{4},1]$ and $T_2 \times [\frac{3}{4},1]$ the unions of
the $m$ basic slices in $T_1 \times [\frac{1}{2},1]$ and $T_2
\times [\frac{1}{2},1]$ closest to $T_1 \times \{1\}$ and $T_2
\times \{1\}$, respectively. Remove from $M$ the open solid tori
bounded by $T_1 \times\{\frac{3}{4}\}$, $T_2
\times\{\frac{3}{4}\}$ and $T_3 \times\{1\}$, and apply Lemma
\ref{GS36} to the resulting $\Sigma \times S^1$ and the thickened
tori $T_1 \times [\frac{3}{4},1]$ and $T_2 \times
[\frac{3}{4},1]$. Then there exists a properly embedded convex
vertical annulus $A$ in $\Sigma \times S^1$ with no
$\partial$-parallel dividing curves, whose Legendrian boundary
$\partial A=(A\cap (T_1 \times\{\frac{3}{4}\}))\cup (A\cap (T_2
\times\{\frac{3}{4}\}))$ intersects the dividing sets of $T_1
\times\{\frac{3}{4}\}$ and $T_2 \times\{\frac{3}{4}\}$
efficiently. Cutting $\Sigma \times S^1$ open along $A$ and round
the edges, we get a thickened torus $T_3 \times [1,2]$ with convex
bounary such that $T_3 \times \{1\}$ has two dividing curves of
slope $\infty$, and $T_3 \times \{2\}$ has two dividing curves of
slope $\frac{1}{m}$. Note that $\frac{1}{m}\leq-\frac{q_3}{p_3}$.
If $\frac{1}{m}=-\frac{q_3}{p_3}$, then, as above, the existence
of $T_3 \times [1,2]$ means that we can thicken $V_3$ to
$\widetilde{V}_3$ such that $\xi|_{\widetilde{V}_3}$ is
overtwisted, which contradicts the tightness of $\xi$. If
$\frac{1}{m}<-\frac{q_3}{p_3}$, we can choose $s_3$ so that
$\frac{1}{m}<s_3<-\frac{q_3}{p_3}$. Then the thickened torus $T_3
\times [0,2]=(T_3 \times I) \cup (T_3 \times [1,2])$ has
$I$-twisting greater than $\pi$, which again contradicts the
tightness of $\xi$. So $T_2 \times [\frac{1}{2},1]$ can not
contain $m$ basic slices of the sign $-\sigma$.

\textit{Case 2.} Assume that $T_2 \times [\frac{1}{2},1]$ contains
$m+1$ basic slices of the sign $\sigma$. We shuffle one of the
$\sigma$ sign to the basic slice in $T_2 \times [\frac{1}{2},1]$
closest to $T_2 \times \{1\}$,and denote by $T_2 \times
[\frac{3}{4},1]$ this basic slice. Similar to the proof of Theorem
\ref{-2}, we can find a convex vertical annulus $A$ in $M$
satisfying:
\begin{enumerate}
    \item $A$ has no $\partial$-parallel dividing curves;
    \item $\partial A=(A\cap (T_2 \times\{\frac{3}{4}\}))\cup
    (A\cap (T_3 \times\{\frac{1}{2}\}))$, which is Legendrian and
    intersects the dividing sets of $T_2 \times\{\frac{3}{4}\}$ and
    $T_3 \times\{\frac{1}{2}\}$ efficiently;
    \item $A$ is disjoint from $T_1$ and the interior of the solid tori in $M$ bounded
    by $T_2 \times\{\frac{3}{4}\}$ and $T_3
    \times\{\frac{1}{2}\}$.
\end{enumerate}

Note that, since $q_1=1$, the diffeomorphism $\varphi_1:\partial
V_1\rightarrow T_1$ is given by
\[
\varphi_1 ~=~
\left(%
\begin{array}{cc}
  p_1 & 1 \\
  -1 & 0 \\
\end{array}%
\right).
\]

Remove from $M$ the interior of the solid tori in $M$ bounded by
$T_2 \times\{\frac{3}{4}\}$ and $T_3 \times\{\frac{1}{2}\}$, and
cut it open along $A$. We get a thickening $\widetilde{V}_1$ of
$V_1$, whose boundary is convex with two dividing curves of slope
$\infty$. Then $\widetilde{V}_1$ is a standard neighborhood of a
Legendrian circle isotopic to the $\frac{1}{p_1}$-singular fiber
with twisting number $0$. We stabilize this Legendrian circle
once. This gives a thickened torus $\widetilde{T}_1\times[0,2]$
with convex boundary such that
$\widetilde{T}_1\times\{2\}=\varphi_1(\partial\widetilde{V}_1)$,
which has two dividing curves of slope $0$, and
$\widetilde{T}_1\times\{0\}$ has two dividing curves of slope
$-\frac{1}{p_1-1}$. Since we can choose the sign of the
stabilization, we can make the sign of this basic slice $\sigma$.
Since $-\frac{1}{p_1-1}\geq-\frac{1}{2m}$, we can find tori
$\widetilde{T}_1\times\{\frac{1}{2}\}$ and
$\widetilde{T}_1\times\{1\}$ in $\widetilde{T}_1\times[0,2]$
parallel to the boundary such that
$\widetilde{T}_1\times\{\frac{1}{2}\}$ has two dividing curves of
slope $-\frac{1}{2m}$, and $\widetilde{T}_1\times\{1\}$ has two
dividing curves of slope $\infty$. Note that
$\widetilde{T}_1\times[\frac{1}{2},1]$ is now a continued fraction
block consisting of $2m$ basic slices of the sign $\sigma$. Now
use $\widetilde{T}_1\times\{1\}$ to thicken $T_2 \times
[\frac{1}{2},\frac{3}{4}]$ to $\widetilde{T}_2 \times
[\frac{1}{2},1]$ such that $\widetilde{T}_2 \times
[\frac{1}{2},\frac{3}{4}]=T_2 \times [\frac{1}{2},\frac{3}{4}]$,
and $\widetilde{T}_2 \times \{1\}$ has two vertical dividing
curves. Note that $\widetilde{T}_2 \times [\frac{1}{2},1]$ is a
continued fraction block that contains at least $m$ basic slices
of the sign $\sigma$. Now, similar to Case 1, we can find a
contradiction. Thus, $T_2 \times [\frac{1}{2},1]$ can not contain
$m+1$ basic slices of the sign $\sigma$ either.

But $T_2 \times [\frac{1}{2},1]$ contains $2m$ basic slices. So
either $m$ of these are of the sign $-\sigma$, or $m+1$ of these
are of the sign $\sigma$. This is a contradiction. Thus, no such
$\xi$ exists on $M$, and, hence, we proved part (3) of Theorem
\ref{-1}. \hfill $\Box$

\vspace{.4cm}

\section{Final Remarks}

When a Seifert space $M$ has more than $3$ singular fibers, or the
base surface has genus $\geq1$, there exists a vertical
incompressible torus $T$ embedded in $M$. Using this vertical
incompressible torus $T$, we can construct infinitely many
pairwise non-isomorphic universally tight contact structure on
$M$. (See \cite{HKM} for more details.) For infinitely many
universally tight contact structures constructed this way, there
is a tubular neighborhood $T\times I$ of $T$ in $M$ that has
$I$-twisting equal to $\pi$. In such a neighborhood, we can always
find a convex torus $T'$ isotopic to $T$ whose dividing set
consists of vertical circles. By the Legendrian Realization
Principle, we can realize a vertical circle on $T'$ disjoint from
the dividing set as a Legendrian curve. This Legendrian curve is a
vertical Legendrian circle with twisting number $0$. So any such
"larger" Seifert space admit tight contact structures for which
there exists a Legendrian vertical circle with twisting number
$0$.

Using Theorem \ref{-2} and Ghiggini and Sch\"{o}nenberger's result
that, when $r\leq\frac{1}{5}$, no tight contact structures on the
small Seifert space $M(r,\frac{1}{3},-\frac{1}{2})$ admit
Legendrian vertical circles with twisting number $0$, one can
easily prove that the Brieskorn homology spheres
$\pm\Sigma(2,3,p)$ do not admit tight contact structures for which
there exists Legendrian vertical circles with twisting number $0$.
It's certainly very interesting to see if this is true for all
small Brieskorn homology spheres. Unfortunately, Theorem \ref{-1}
is not strong enough to apply to this case. We will have to
develop new techniques before we can tackles this problem.

\end{document}